\pdfoutput=1
\RequirePackage{ifpdf}
\ifpdf 
\documentclass[pdftex]{sigma}
\else
\documentclass{sigma}
\fi

\usepackage{tikz-cd}
\usepackage{setspace}

\numberwithin{equation}{section}

\newtheorem{Theorem}{Theorem}[section]
\newtheorem*{Theorem*}{Theorem}

 { \theoremstyle{definition}
\newtheorem{Definition}[Theorem]{Definition}

 }

\begin{document}

\renewcommand{\thefootnote}{}

\newcommand{\arXivNumber}{2308.02609}

\renewcommand{\PaperNumber}{045}

\FirstPageHeading

\ShortArticleName{The Cobb--Douglas Production Function and the Old Bowley's Law}

\ArticleName{The Cobb--Douglas Production Function\\ and the Old Bowley's Law\footnote{This paper is a~contribution to the Special Issue on Symmetry, Invariants, and their Applications in honor of Peter J.~Olver. The~full collection is available at \href{https://www.emis.de/journals/SIGMA/Olver.html}{https://www.emis.de/journals/SIGMA/Olver.html}}}

\Author{Roman G.~SMIRNOV~$^{\rm a}$ and Kunpeng WANG~$^{\rm b}$}

\AuthorNameForHeading{R.G.~Smirnov and K.~Wang}

\Address{$^{\rm a)}$~Department of Mathematics and Statistics, Dalhousie University, \\
\hphantom{$^{\rm a)}$}~6297 Castine Way, PO BOX 15000, Halifax, Nova Scotia, B3H 4R2, Canada}
\EmailD{\href{mailto:roman.smirnov@dal.ca}{roman.smirnov@dal.ca}}

\Address{$^{\rm b)}$~Sichuan University-Pittsburgh Institute (SCUPI), Sichuan University,\\
\hphantom{$^{\rm b)}$}~610207 Chengdu, Sichuan, P.R.~China}
\EmailD{\href{mailto:kunpeng.wang@scupi.cn}{kunpeng.wang@scupi.cn}}

\ArticleDates{Received July 31, 2023, in final form May 13, 2024; Published online May 30, 2024}

\Abstract{Bowley's law, also referred to as the law of the constant wage share, was a~noteworthy empirical finding in economics, suggesting that a nation's wage share tended to remain stable over time, as observed through most of the 20th century. The wage share represents the proportion of a country's economic output that is distributed to employees as compensation for their labor, usually in the form of wages. The term ``Bowley's law'' was coined in 1964 by Paul Samuelson, the first American laureate of the Nobel memorial prize in economic sciences. He attributed this principle to Sir Arthur Bowley, an English economist, mathematician, and statistician. In this paper, we introduce a mathematical model derived from data for the American economy, originally employed by Cobb and Douglas in 1928 to validate the renowned Cobb--Douglas production function. We utilize symmetry methods, particularly those developed by Peter Olver, to elucidate the validity of Bowley's law within our model's framework. By employing these advanced mathematical techniques, our objective is to elucidate the factors contributing to the stability of the wage share over time. We demonstrate that the validity of both Bowley's law and the Cobb--Douglas production function arises from the robust growth of an economy, characterized by expansion in capital, labor, and production, which can be approximated by an exponential function. Through our analysis, we aim to offer valuable insights into the underlying mechanisms that support Bowley's law and its implications for comprehending income distribution patterns in economies.\looseness=-1}


\Classification{34C14; 37C79; 62P20; 91B39; 91B55; 91B62}

\begin{flushright}
{\em ``Nature imitates mathematics.''}
--- Gian-Carlo Rota
\end{flushright}

\renewcommand{\thefootnote}{\arabic{footnote}}
\setcounter{footnote}{0}

\section{Introduction}
\label{s1}

Throughout his illustrious and productive career, which spans several decades and is still in full swing, Peter Olver has made numerous contributions to the development of the theories of Lie groups, symmetries of differential equations, invariants, and their applications. His work has had a far-reaching impact, spanning fields from image recognition to mathematical physics and geometry (see, for example, \cite{Olver6, Olver6a, Olver4, Olver11, Olver13, Olver1, Olver2, Olver3, Olver7, Olver8, Olver9, Olver10, Olver12,Olver14, Olver15, Olver5}). Additionally, he has been an inspiring, supportive, and mentoring figure, guiding hundreds of mathematicians at various stages of their careers.

We have been privileged to learn from Peter through his lectures, numerous books and articles, which have been particularly instrumental in the first author's (RGS) contributions to the theories of bi-Hamiltonian systems and the invariant theory of Killing tensors (refer to~\cite{CMS17,HMS05,PS05} and the relevant references therein). The latter has been significantly influenced by the monograph Olver \cite{Olver7}.

Our main goal is to demonstrate that Peter Olver's approach, which encompasses symmetry methods, invariant theory, and the theory of Lie groups for studying differential equations, can be readily applied to the investigation of economic models. By demonstrating this compatibility, it opens up new avenues for using advanced mathematical tools in the field of economics. Such interdisciplinary connections have the potential to enhance our understanding of economic phenomena and improve the accuracy of economic models.

Traditionally, when it comes to applications, symmetry methods have demonstrated their value in the study of differential equations, typically derived from Newton's second law or the Lagrangian (Hamiltonian) of the system in question. In this context, a differential equation or a system of differential equations to be studied is derived from first principles using physics.

On the other hand, in the field of economics, among many others, this approach cannot always be applied to derive differential equations that can be employed to model economic phenomena. Instead, differential equations in this setting are often derived from available data rather than first principles. However, we firmly believe that this should not be seen as an obstacle to utilizing symmetry methods and related techniques to study economic models, particularly those where conserved quantities play a central role. This approach can be a powerful tool, enabling us to identify and analyze key features that contribute to a better understanding of conserved quantities in economic systems. In turn, it could provide valuable insights into the dynamics and stability of such models, leading to a deeper comprehension of economic behaviors and potentially facilitating better decision-making in economic policy and planning.

It must be noted that significant efforts have already been made to incorporate symmetry methods and the methods of Lie group theory into the field of mathematical modeling in economics. First and foremost, we wish to mention the influential books by Sato \cite{Sato81} and Sato and Ramachandran \cite{SR14} that have contributed significantly to this endeavor. Additionally, relevant contributions can be found, for example, in the works of Perets and Yashiv \cite{PY18}, Fukang \cite{Fukang96}, and others \cite{SW20}. This comes as no surprise because many basic semi-heuristic economic laws can be derived, for example, from the corresponding scale or shift-invariance under an appropriate symmetry transformation \cite{KKTT19}.

In the following discussion, we will revisit the so-called ``Bowley's law'' \cite{Bowley1900, Bowley1937, PS1964}. Paul Samuelson introduced the term ``Bowley's law" in 1964, within the sixth American edition of his renowned textbook {\em Economics} \cite{PS1964}, to describe the constant wage share observed in economic data. The term was a tribute to Sir Arthur Bowley, an English economist, mathematician, and statistician credited with pioneering the systematic collection and statistical analysis of wage data in the UK. Bowley had hypothesized as early as 1920 that the wage share might remain constant, and later, with Josiah Stamp, provided evidence supporting this hypothesis by comparing UK wage shares in 1911 and 1924. In his 1937 book {\em Wages and Income in the United Kingdom since $1860$} \cite{Bowley1937}, Bowley definitively stated the constancy of the wage share, a~discovery that challenged the views of classical economists such as Ricardo. They traditionally regarded the proportions of land, capital, and labor as intrinsically variable.

This observation made by Arthur Bowley has been viewed as a purely heuristic insight that has been the subject of much controversy since its inception to the present day. The economic data collected in various countries from the end of the 19th century until about 1980 gave rise to and strongly supported this law, leading to its wide acceptance within the economics community at the time. However, recent data patterns appear to deviate from this law, casting doubt on its continued validity (see Kr\"{a}mer \cite{Kramer11} for more details and references). Our objective is to build upon the work initiated in \cite{SW20} to construct a well-founded mathematical model. This model aims to shed light on the validity (or lack thereof) of Bowley's law in economics and explore how it manifests itself, starting from the available economic data and extending all the way to the corresponding fundamental invariants of the appropriate one-parameter Lie group transformations. These transformations will be obtained in the form of a data-driven dynamical system extracted from the said economic data.

This paper is structured as follows. In Section \ref{s2}, we provide a mathematical definition of wage share by deriving the corresponding formula from an optimization problem. Section \ref{s3} establishes a connection between the Cobb--Douglas production function and Bowley's law. In Section \ref{s4}, we derive the Cobb--Douglas function from an exponential model using mathematical and statistical methods. In Section \ref{s5}, we demonstrate the notion of labor (wage) share as a~time-independent invariant of a prolonged infinitesimal one-parameter Lie group action, utilizing symmetry methods.
In Section \ref{s6}, we introduce the logistic model as a generalization of the exponential model. The material for Section \ref{s7} demonstrated that Bowley's law no longer valid within the framework of the logistic model. Finally, in Section \ref{s8}, we present our concluding remarks.\looseness=-1

\section{What is wage share?}
\label{s2}

Our model will revolve around the notion of a production function. In economics, a production function represents the relationship between physical output and input factors, which have historically been considered by mainstream economists to include capital, labor, land, and entrepreneurship. In most modern economic models, a production function is a function of two input factors: capital and labor. However, other schools of thought propose alternative perspectives. For instance, some consider a production function to be a function of energy generated within the framework of a given economy (see, for example, \cite{CS21}), while others link it to the maximum machine speeds \cite{Beaudreau17}.

In light of the above, in the subsequent analysis, we will define our configuration space using the following three variables: production ($Y$), labor ($L$), and capital ($K$). Furthermore, we assume that $Y$ (output) is a function of both $L$ and $K$ (input factors):
\begin{equation}
\label{1}
Y = f(L, K).
\end{equation}

Next, we recall first that wage share, also known as labor share, refers to the portion of national income that is allocated to workers as compensation for their labor. It represents the percentage or fraction of total income in an economy that goes to wages, salaries, and other forms of remuneration received by workers for their work.

The wage share is calculated by dividing the total wages and compensation paid to workers by the total output or gross domestic product (GDP) of the economy. It is an essential indicator used to understand the distribution of income between labor and capital in an economy. A~high wage share suggests that a~larger portion of national income is going to workers, while a~low wage share indicates a larger share is going to other factors of production, such as capital or profits. Changes in the wage share over time can reflect shifts in the dynamics of income distribution and have implications for economic inequality and overall economic well-being.

To obtain a formula for the wage share in an economy governed by production function \eqref{1} that is assumed to operate ideally, let us consider the following unconstrained optimization problem:
\begin{equation}
\label{2}
\Pi = p Y - wL - rK \rightarrow \max,
\end{equation}
where $\Pi$ is the profit with nominal wage ($w$), nominal rent ($r$), and nominal price ($p$). We assume $Y = f(K, L)$ is continuously differentiable and such that there is an interior solution for~${K, L, Y \ge 0}$ and so, in particular, we have
\begin{equation}\label{3}
\frac{\partial Y}{\partial L} = \frac{w}{p}.
\end{equation}
Then, in view of the above definition and \eqref{3}, wage (labor) share is easily found to be
\begin{equation}\label{4}
s_L = \frac{wL}{pY} = \frac{\partial Y}{\partial L}\frac{L}{Y}.
\end{equation}
Similarly, we conclude that the corresponding formula for capital share $s_K$ is given by
\begin{equation}\label{5}
s_K = \frac{rK}{pY} = \frac{\partial Y}{\partial K}\frac{K}{Y}.
\end{equation}

We will establish, with the necessary mathematical rigor and employing symmetry methods, that under realistic assumptions regarding the evolution of the three variables $Y$, $L$, and $K$, the quantities $s_L$ and $s_K$ are constant over time, serving as time-independent invariants. This constitutes the essence of Bowley's law, asserting that the wage share remains unchanged over time. Until approximately 1980, economic data strongly supported this law, and it was widely embraced by the economics community throughout most of the 20th century, starting from its introduction by Arthur Bowley. However, in recent times, the support for this law has waned, as evident in sources like Kr\"{a}mer \cite{Kramer11}, Orlando \cite{OG23}, Schneider \cite{Schneider11}, and Stockhammer \cite{SE13}, which provide further details and references on this topic.

\section{Bowley's law and the Cobb--Douglas production function}
\label{s3}

It is well known that Bowley's law holds true when the production function \eqref{1} used in formula~\eqref{2} corresponds to the celebrated Cobb--Douglas production function \cite{CD28}, which is given~by%
\begin{equation}
\label{6}
Y = AL^{\alpha}K^{\beta},
\end{equation}
where $Y = f(L, K)$, $L$, and $K$ as defined before, $A$ is total factor productivity, while $\alpha$ and $\beta$ are the output elasticities of capital and labor respectively. In 1928, Charles Cobb and Paul Douglas presented a paper~\cite{CD28} (see also \cite{Douglas76}) focused on studying the growth of the American economy between 1899 and 1922. To model the production output, they employed the function~\eqref{6} previously introduced and studied by Knut Wicksell, Philip Wicksteed, and L\'{e}on Walras (see Humphrey \cite{Humphrey97}. Of particular importance for mathematical modeling in economics is the case when the Cobb--Douglas function \eqref{6} has constant returns to scale, namely when
\begin{equation}
\label{7}
\alpha+\beta = 1.
\end{equation}
In mathematical terms, this signifies that $Y$ as a function of $L$ and $K$ is a homogeneous function of degree one. From an economic perspective, this implies that when both capital and labor in \eqref{6} increase by the same factor $\lambda$, the corresponding increase in $Y$ will also be $\lambda$. In other words, if both inputs are scaled up proportionally, the output will increase by the same proportion. Through the method of least squares estimation, Cobb and Douglas in \cite{CD28} assumed constant returns to scale in \eqref{6}, meaning
\[Y = A L^{\alpha}K^{1-\alpha},\]
 and obtained a labor exponent result of $\alpha = 0.75$, which was later verified by the National Bureau of Economic Research to be~$0.741$. In their subsequent work in the 1940s, they allowed for variable exponents on $K$ and $L$, leading to estimates that closely matched improved measures of productivity developed during that period (see for more details). The constant $A$ (total productivity) was determined to be $A = 1.01$. In the final analysis, after observing how the Cobb--Douglas function was derived from a data set \cite{CD28}, it becomes apparent that the process was rooted in data mining. As a result, we can succinctly summarize it as follows:
$$
\fbox{\rule[-3mm]{0cm}{0.8cm}Data} \rightarrow \fbox{\rule[-3mm]{0cm}{0.8cm}Production function}
$$

At this point one might naturally ask a direct question: Are there any other functions in the form of \eqref{6} that offer equally accurate fits to the data representing the growth of the American economy from 1899 to 1922, comparable to the function used by Cobb and Douglas in \cite{CD28}? The answer to this question is ``Yes'', and we will soon elaborate on the meaning and importance of this observation.

Under the assumption that the production function $Y$ in \eqref{2} is in the form \eqref{6}, it is easy to verify the validity of Bowley's law. Indeed, substituting \eqref{6} in \eqref{4}, it follows that in this case $s_L = \alpha$,
that is wage (labor) share is constant and equal to the output elasticity of labor~$\beta$. Similarly, we get for capital share $s_K = \beta$.

Therefore, based on our analysis, we can conclude that the Cobb--Douglas function, describing the growth of production according to the formula \eqref{6}, implies Bowley's law.

\section[The Cobb--Douglas function derived from an exponential model]{The Cobb--Douglas function derived\\ from an exponential model}
\label{s4}

Based on the information presented in Section \ref{s3}, it can be concluded that Cobb and Douglas established in
\cite{CD28} the legitimacy of the production function \eqref{6} through a combination of economic reasoning, statistical methods, and empirical data. Specifically, they demonstrated that by assigning the values $\alpha = 0.75$, $\beta = 0.25$, and $A = 1.01$ to the respective parameters, the production function enjoying constant returns to scale yielded a precise fit to the data set under consideration.

Sato \cite{Sato81} pursued a similar objective, combining economic reasoning with analytical methods to directly derive the Cobb--Douglas function without any reliance on statistical data. Specifically, Sato \cite{Sato81} (also referenced in \cite{SW20}) employed Lie group theory methods to establish the Cobb--Douglas function as an invariant of a two-dimensional integrable distribution, determined by exponential growth in labor, capital, and production. In the subsequent discussion, we shall integrate both approaches, making necessary simplifications in the process.

Specifically, in \cite{Sato81}, Sato introduced the concept of simultaneous holotheticity, which describes a situation where two sectors of the same economy are governed by the same aggregate production function, while experiencing different levels of technical change within each sector. Assuming that the growth in $K$ (capital), $L$ (labor), and $Y$ (production) is exponential in both sectors and working at the level of infinitesimal action, this assumption leads to the following distribution of vector fields for which we seek an invariant function $\varphi (K, L, Y)$:
\begin{gather}
X_1\varphi =b_1K \frac{\partial \varphi}{\partial K} + b_2 L \frac{\partial\varphi }{\partial L} + b_3Y\frac{\partial \varphi}{\partial f}=0, \nonumber \\
X_2\varphi =b_4K \frac{\partial \varphi}{\partial K} + b_5 L \frac{\partial\varphi }{\partial L} + b_6Y\frac{\partial \varphi}{\partial f}=0.
\label{Hol}
\end{gather}
To ensure that the resulting invariant Cobb--Douglas function \eqref{6} satisfies equation \eqref{7}, Sato assumed the following values for the parameters $b_i$ where $i = 1, \ldots, 6$: $b_1 = b_2 = b_3 = b_6 = 1$, $b_4 = a > 0$, and $b_5 = b > 0$. It is evident that the vector fields $X_1$ and $X_2$ form a two-dimensional integrable distribution on $\mathbb{R}_+^2$. By integrating the corresponding total differential equation,
\[(YL - bYL){\rm d}K + (aYK-YK){\rm d}L + (bKL - aKL){\rm d}Y = 0,\]
or
\begin{gather*}
(1-b)\frac{{\rm d}K}{K} + (a-1)\frac{{\rm d}L}{L} + (b-a)\frac{{\rm d}Y}{Y} = 0,
\end{gather*}
one arrives at the Cobb--Douglas function given by
\begin{equation}
\label{CD}
Y = AL^{\alpha} K^{\beta},
\end{equation}
where the output elasticities
\begin{equation}
\label{oe}
\alpha = \frac{1-b}{a-b}, \qquad \beta = \frac{a-1}{a-b}
\end{equation}
satisfy the condition of constant returns to scale \eqref{7}, as desired. However, we note that the condition of simultaneously holotheticity is rather artificial.

Moreover, it is worth noting that the output elasticities $\alpha$ and $\beta$ given by \eqref{oe} are strictly positive if and only if either $a > b$, $b < 1$, and $a > 1$, or $a < b$, $b > 1$, and $a < 1$. This implies that assuming exponential growth in $K$ (capital), $L$ (labor), and $Y$ (production) leads to the derivation of the Cobb--Douglas function. In other words, the Cobb--Douglas function arises as a~consequence of exponential growth. However, we argue that working directly with the Lie group action (as seen in Olver \cite{Olver7}) is much more straightforward than dealing with a distribution of vector fields and subsequently integrating it. Additionally, we believe that to derive the Cobb--Douglas function \eqref{6} with constant returns to scale, there is no need to assume simultaneous holotheticity.

Furthermore, a crucial question remains: ``Are the approaches to the derivation of the production function \eqref{6}, as proposed by Cobb--Douglas and Sato, compatible?''

Following Sato, let us assume that labor, capital, and production grow exponentially and can be treated as functions of time. This assumption leads to the following very simple dynamical system:
\begin{equation}
\label{model1}
\dot{x}_i=b_ix_i, \qquad i = 1,2,3,
\end{equation}
where $x_1 = x_1(t) = L(t)$ (labor), $x_2 = x_2(t) =K(t)$ (capital), $x_3= x_3(t) = Y(t)$ (production).

Noticing that the system \eqref{model1} can be viewed as a one-parameter Lie group action in $\mathbb{R}^3_+$, we get from \eqref{model1}
\begin{equation}
\label{action}
x_i = x^0_i {\rm e}^{b_it}, \qquad i = 1, 2, 3,
\end{equation}
where the parameters $b_i$ are as before and $x_i^0$ are the constants of integration (initial conditions). Now, let us aim to eliminate the parameter $t$. Consider the following product:
\begin{equation}
\label{inv1}
\big(x^0_1{\rm e}^{b_1t}\big)^{a_1} \big(x^0_2{\rm e}^{b_2t}\big)^{a_2} \big(x^0_3{\rm e}^{b_3t}\big)^{a_3} =\big(x^0_1\big)^{a_1}\big(x^0_2\big)^{a_2}\big(x^0_3\big)^{a_3} {\rm e}^{(a_1b_1 + a_2 b_2 + a_3 b_3)t}.
\end{equation}
It is easy to see that the expression \eqref{inv1} is independent of $t$ and thus defines an invariant of the group action \eqref{action}, provided that the following ``orthogonality'' condition holds true:
\begin{equation}
\label{linear}
a_1 b_1 + a_2 b_2 + a_3 b_3 = 0 \ \Leftrightarrow \ {\bf a} \cdot {\bf b} = 0,
\end{equation}
where ${\bf a} = \langle a_1, a_2, a_3\rangle$, ${\bf b} = \langle b_1, b_2, b_3\rangle$. Note the condition \eqref{linear} defines a plane in $\mathbb{R}^3$ by the normal vector ${\bf b}$. Therefore, it follows from \eqref{inv1} that the function
\begin{equation}
\label{inv2}
f(x_1, x_2, x_3) = \big(x^0_1\big)^{a_1}\big(x^0_2\big)^{a_2}\big(x^0_3\big)^{a_3} x_1^{a_1}x_2^{a_2} x_3^{a_3}
\end{equation}
is constant along the flow generated by \eqref{action} if and only if the condition \eqref{linear} is satisfied. Hence, we conclude that the orthogonality condition \eqref{linear} is the true invariant that determines the family of Cobb--Douglas production functions defined by the formula \eqref{inv2}. Indeed, we have utilized the method of moving frames \cite{Olver6, Olver6a, Olver7}, albeit at a rudimentary level.

 Consider now a particular level set
\begin{equation}
\label{inv3}
\big(x^0_1\big)^{a_1}\big(x^0_2\big)^{a_2}\big(x^0_3\big)^{a_3} x_1^{a_1}x_2^{a_2} x_3^{a_3} = C,
\end{equation}
for a fixed constant $C \not=0$. Solving \eqref{inv3} for $x_3$ (production), we arrive at
\begin{equation}
\label{CDF1}
x_3 = f(x_1, x_2) = \left(\frac{C}{\big(x^0_1\big)^{a_1}\big(x^0_2\big)^{a_2}\big(x^0_3\big)^{a_3}}\right)^{\frac{1}{a_3}}x_1^{-\frac{a_1}{a_3}}x_2^{-\frac{a_2}{a_3}}.
\end{equation}
Setting
\[
A = \left(\frac{C}{\big(x^0_1\big)^{a_1}\big(x^0_2\big)^{a_2}\big(x^0_3\big)^{a_3}}\right)^{\frac{1}{a_3}},\qquad {\alpha} = -\frac{a_1}{a_3},\qquad {\beta} = -\frac{a_2}{a_3},
\]
 identifying $x_3 = Y$ (production), $x_1 = L$ (labor), $x_2 = K$ (capital), we arrive at the Cobb--Douglas production function \eqref{6}, provided (and that is the key!)
\begin{equation}
\label{linear1}
\alpha b_1 + \beta b_2 - b_3 = 0,
\end{equation}
where $b_i$, $i = 1, 2, 3$ are determined by the exponential growth in the input factors and production given by \eqref{action}. The equation \eqref{linear1} defines a straight line in the $\alpha\beta$-plane that intersects the line \eqref{7} in the first quadrant, provided $\alpha, \beta, b_1, b_2, b_3 >0$ and either of the two inequalities
\begin{gather}
\label{ineq}
b_2>b_3>b_1,\\
\label{ineq1}
b_1>b_3>b_2,
\end{gather}
hold true. In fact, it follows from \eqref{6} and \eqref{linear1} that we have a {\em family of production functions of the Cobb--Douglas type} determined by (see also~\cite{SWW22})
\begin{equation}
\label{family}
Y = f(L, K) = A L^{\alpha} K^{\frac{b_3}{b_2} - \alpha\frac{b_1}{b_2}}, \qquad \alpha>0.
\end{equation}
We conclude, therefore, that the one-parameter Lie group action \eqref{action} admits a family of the Cobb--Douglas functions given by \eqref{family} if and only if the output elasticities $\alpha$ and $\beta$ are constrained by the orthogonality condition \eqref{linear1}. Moreover, if additionally the parameters~$b_i$,~${i=1,2,3}$ satisfy either the inequality \eqref{ineq} or \eqref{ineq1}, we can always use the formula~\eqref{linear1} to pick among the functions given by \eqref{6}, the Cobb--Douglas function enjoying constant
returns to scale \eqref{7}.

\begin{figure}[ht]\centering
\includegraphics[width=0.5\textheight]{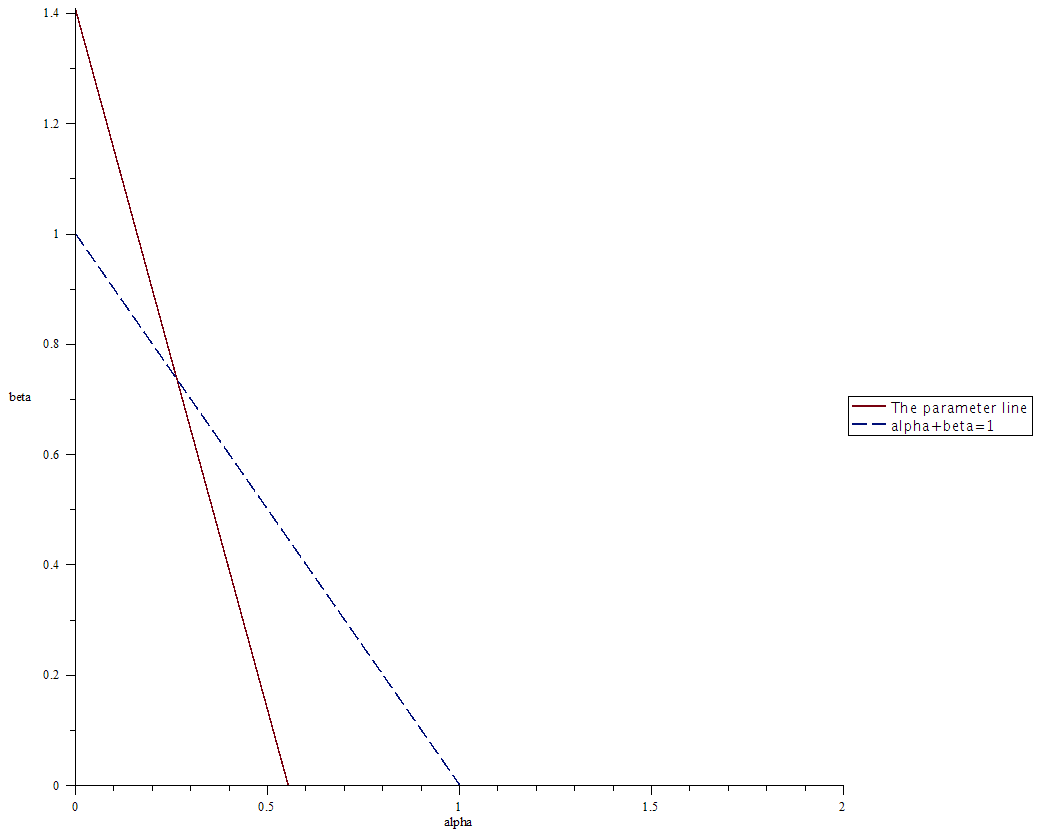}
\caption{Constant returns to scale is possible.}\label{plot1}
\end{figure}

\begin{figure}[ht]\centering
\includegraphics[width=0.5\textheight]{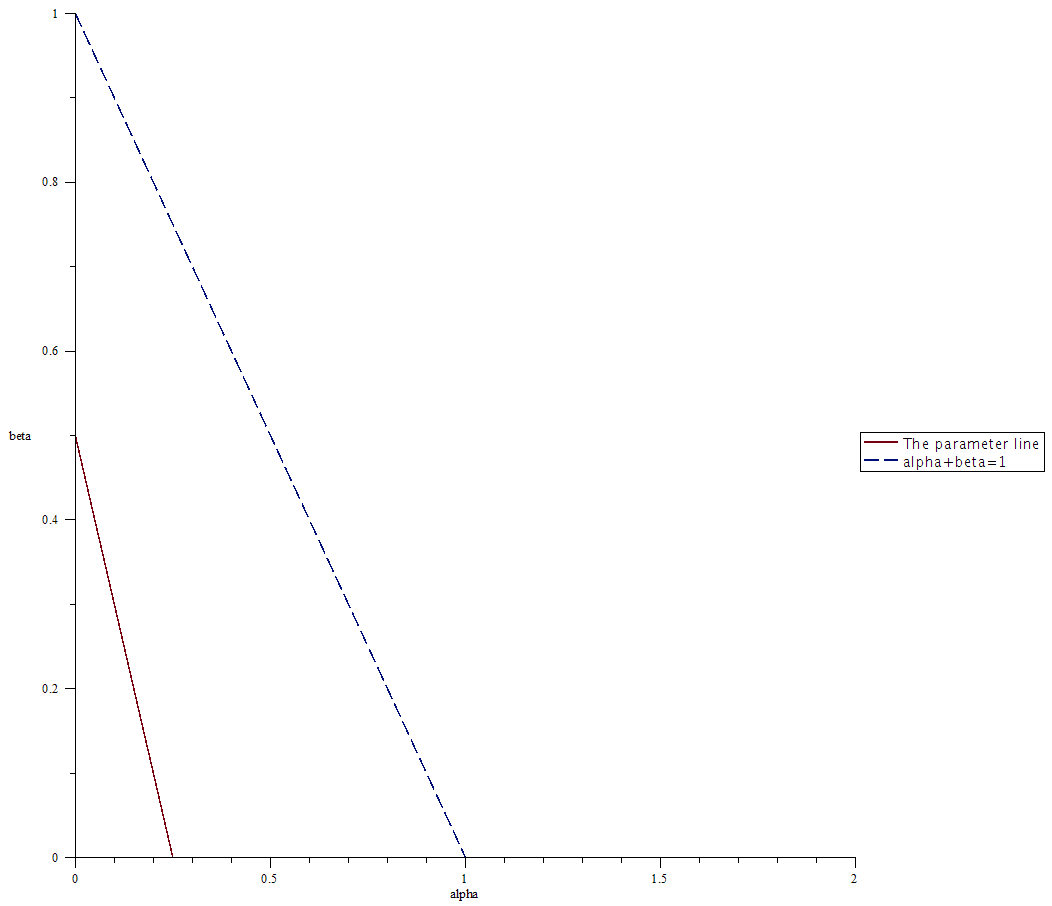}
\caption{Decreasing return to scale ($\alpha + \beta <1$).}\label{plot2}
\end{figure}

\begin{figure}[ht]\centering
\includegraphics[width=0.5\textheight]{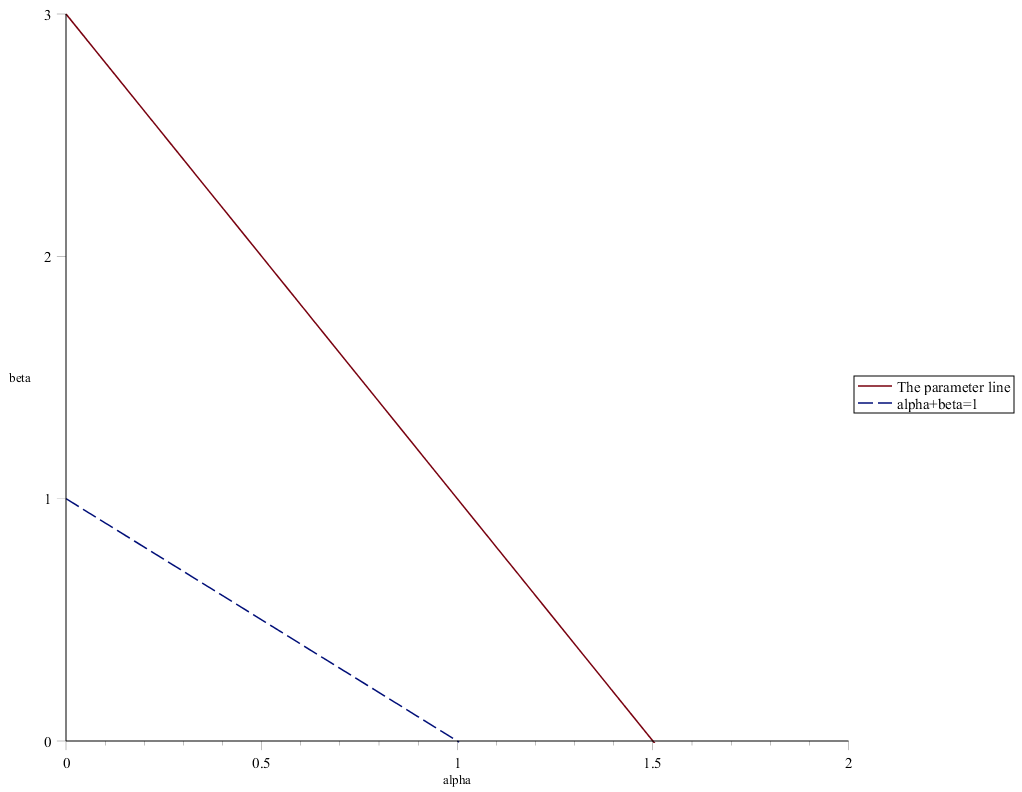}
\caption{Increasing return to scale ($\alpha + \beta >1$).}\label{plot3}
\end{figure}

In summary, we have the following three cases determined by the linear equations \eqref{7} and~\eqref{linear1} and illustrated by Figures \ref{plot1}, \ref{plot2}, and \ref{plot3}, respectively.

We observe that Sato \cite{Sato81} chose the values of the parameters $b_i$, $i = 1, \ldots, 6$ in \eqref{Hol} to assure that the two lines defined by \eqref{linear1} that corresponded to the vector fields $X_1$ and $X_2$ given by~\eqref{Hol} necessarily intersected in the first quadrant of the $\alpha\beta$-plane (see Figure \ref{plot1}). Indeed, it is easy to see that the lines $\alpha + \beta - 1 =0$ (the vector field $X_1$) and $a\alpha +b\beta - 1 = 0$ (the vector field $X_2$) intersect in the first quadrant of the $\alpha\beta$-plane, provided the output elasticities \eqref{oe} are strictly positive.

In view of the above observations and results, we propose the following generalization of the Cobb--Douglas function.
\begin{Definition}[Cobb--Douglas function] \label{d1} Given the one-parameter group action
\begin{equation}
\label{action5}
x_i = x^0_i {\rm e}^{b_it}, \qquad x_i^0, b_i>0, \qquad i = 1, \ldots, n
\end{equation}
in $\mathbb{R}^n_+$. Then the {\em Cobb--Douglas function} is defined as an element of the following family of invariants of the one-parameter Lie group action \eqref{action5}:
\begin{equation}
\label{GCDF}
\prod_{i=1}^n \big(x_i^0 x_i\big)^{a_i} = C, \qquad a_i > 0, \qquad i = 1, \ldots, n,
\end{equation}
where $C \in \mathbb{R}$ is an arbitrary constant and $x_i^0$, $i = 1, \ldots, n$ are the corresponding initial conditions, if the orthogonality condition~\eqref{linear6}
holds true.
\end{Definition}

\begin{figure}[t]
\centering
\includegraphics[width=0.5\textheight]{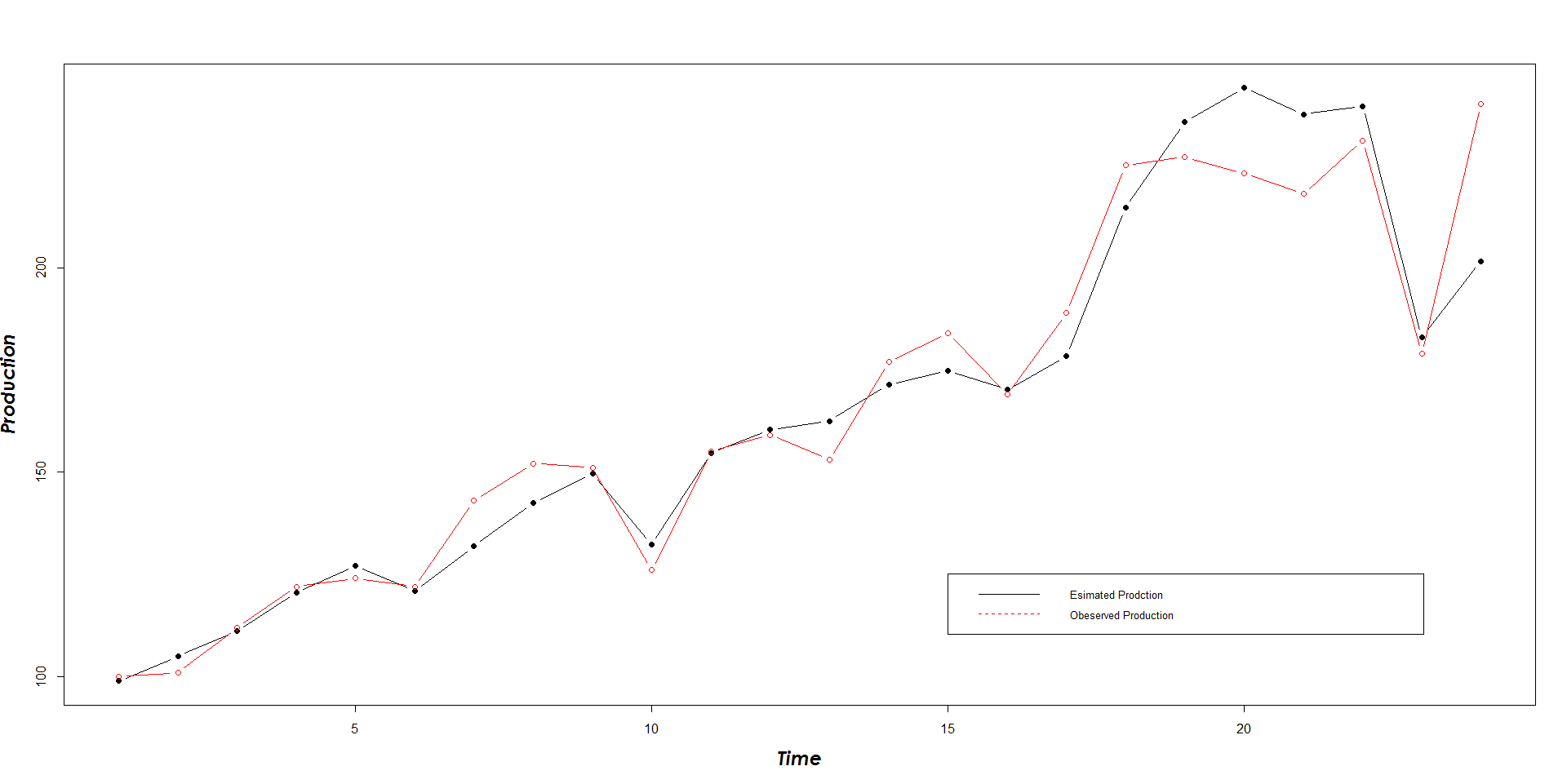}
\caption{The function \eqref{CDF4} vs the index values for the production studied by Cobb and Douglas~\cite{CD28}.}\label{figure1}
\end{figure}

Note that the data examined by Cobb and Douglas in \cite{CD28} corresponds to the case presented in Figure \ref{plot1}. Now, let us review this data in light of the derived formulas and observations made. First, we linearize the equations \eqref{action} by taking the logarithm of both sides of each equation, leading to the following linear expressions:
\begin{equation}
\label{action3}
\ln x_i = C_i + b_i t, \qquad i = 1, 2, 3,
\end{equation}
where $C_i = \ln x^0_i$ for $i = 1, 2, 3$. The authors used the R programming language to fit the formulas \eqref{action} to the index numbers of the industrial output $Y$, fixed capital $K$, and total number of manual workers $L$ on a logarithmic scale, using the data studied by Cobb and Douglas in \cite{CD28} (refer to \cite[Table 1]{SW21}). More specifically, the corresponding coefficients $C_i$ and $b_i$ for~${i=1,2,3}$ in~\eqref{action3} were accurately determined from the data using R and the least squares method, resulting in the following values:
\begin{alignat}{4}
&b_1=0.02549605, \qquad&& C_1=4.66953290 \qquad&&\text{(labor)},& \nonumber\\
&b_2=0.06472564, \qquad&& C_2=4.61213588 \qquad&& \text{(capital)},&\nonumber\\
&b_3=0.03592651, \qquad&& C_3=4.66415363 \qquad&& \text{(production)}.&\label{data}
\end{alignat}
Therefore, it follows that the production, labor and capital represented by the data in \cite{CD28} enjoyed (nearly perfect) exponential growth. Moreover, the coefficients $b_1$, $b_2$, and $b_3$ given by \eqref{data} satisfy the inequality \eqref{ineq}. Thus, we compute the corresponding values of the output elasticities~$\alpha$ and $\beta$ in \eqref{6}, satisfying the condition \eqref{7}, using the formula \eqref{linear1} for the estimated values $b_i$, $i=1,2,3$ given by \eqref{data}. Solving the linear system given by \eqref{7} and \eqref{linear1}, we find%
\begin{equation}
\label{alphabeta}
\alpha = \frac{b_3 - b_2}{b_1 - b_2}, \qquad \beta =\frac{b_3-b_1}{b_2 - b_1}.
\end{equation}
Note, the formula \eqref{alphabeta} was originally derived in \cite{SW21} with the aid of the bi-Hamiltonian approach \cite{Olver4, Olver2, PS05}. Next, substituting the values \eqref{data} into \eqref{alphabeta}, we get
\begin{equation}
\label{alphabeta-1}
\alpha=0.7341175376, \qquad \beta=0.2658824627.
\end{equation}
We see that the values for the output elasticities $\alpha$ and $\beta$ \eqref{alphabeta-1} compare well to the corresponding values found by Cobb and Douglas in \cite{CD28}. However, there are other values of $\alpha$ and $\beta$, satisfying the linear relation \eqref{linear1} for the values of $b_1$, $b_2$, and $b_3$ given by \eqref{data}. For example, setting~${\alpha = 1}$, we find, via \eqref{linear1} and using the values given by \eqref{data}, the corresponding value for $\beta$: \smash{$\beta = \frac{b_3 - b_1}{b_2} = 0.16114881212$}.
Note that the values $\alpha=1$ and $\beta= 0.16114881212$ in this case no longer add up to one, while the function
\begin{equation}
\label{CDF4}
Y = ALK^{0.16114881212}
\end{equation}
is a legitimate Cobb--Douglas function compatible with the data studied in \cite{CD28}, which we have confirmed with the aid of the R programming language. More specifically, we have determined that the value of $A \approx 0.4710156$ affords a good fit of the function \eqref{CDF4} to the data studied by Cobb and Douglas in \cite{CD28}. We compare the time series representing the actual values of the production in the US manufacturing sector for 1899--1922 (see \cite[Table~1]{SW21}) with the values computed according to the formula \eqref{CDF4} in Figure~\ref{figure1} (see also \cite[Chart~II, p.~153]{CD28}). We see that the production function given by \eqref{CDF4} also provides us with a near perfect fit to the data studied by Cobb and Douglas in \cite{CD28}.
We wish to add that the time series given by Figure~\ref{figure1} accurately describes the fluctuations of the US economy of that period. More specifically, we clearly see on it the Panics of 1907 and 1910--11, the Recession of 1913--14, the Post-World War~I Recession of August 1918 -- March 1920, and, finally, the Depression of 1920--21.

Recall that in \cite{SW20} we studied the data from the period 1947--2016 provided by the Federal Reserve Bank of St.\ Louis (\url{https://fred.stlouisfed.org}), employing the FRED tool. We used R to fit the Cobb--Douglas function with constant returns to scale of the form \eqref{6} to the aforementioned data. The best estimation of the Cobb--Douglas function that we managed to have obtained, is as follows:
\begin{equation}
\label{CDF5}
Y= f(K, L) = 0.2464455K^{1.6612365}L^{-0.6612365}.
\end{equation}
Note that although the elasticities of substitution $\alpha$ and $\beta$ add up to one in this case, the parameter $\beta$ is negative and so the form \eqref{CDF5} is incompatible with the definition of the Cobb--Douglas function with constant returns to scale. Next, we use the formulas \eqref{action3} and R to find the corresponding parameters $b_i$, $C_i$, $i =1, 2, 3$:
\begin{alignat}{4}
& b_1=0.06983731,\qquad&& C_1=0.45741448 \qquad && \text{(labor)}, &\nonumber\\
& b_2=0.065705809, \qquad&& C_2=0.75835155 \qquad&& \text{(capital)},&\nonumber\\
& b_3=0.03421333, \qquad&& C_3=2.58402362 \qquad&& \text{(production)}.&\label{data1}
\end{alignat}
We immediately observe that the parameters $b_1$, $b_2$, and $b_3$ in \eqref{data1} do not satisfy either of the inequalities~\eqref{ineq} and~\eqref{ineq1}, which is the reason we cannot employ in this case the formulas~\eqref{alphabeta} to compute the output elasticities~$\alpha$ and $\beta$ satisfying the condition \eqref{7} that defines the Cobb--Douglas function with constant returns to scale~\eqref{7}. Relaxing the constant returns to scale condition~\eqref{7}, we arrive at the following Cobb--Douglas function that fits quite well to the data:
\[
Y = f(K, L) = 9.89921606 L^{0.05018686}K^{0.45529695}.
\]
Indeed, the residual sum of squares (RSS) in this case is 584.4616. Note that the values of $\alpha$ and $\beta$ computed with the aid of the formulas \eqref{alphabeta} and the parameters $b_1$, $b_2$, and $b_3$ in \eqref{data1} approximately satisfy the orthogonality condition \eqref{linear1}.

In summary, for a given dataset, we first utilized statistical tools to derive the corresponding dynamical system that is compatible with the data. Subsequently, employing the appropriate mathematical machinery, we obtained a (family of) production function(s) as time-independent invariants of the said dynamical system. Ultimately, the new scheme can be depicted as follows:
$$
\fbox{\rule[-3mm]{0cm}{0.8cm}Data} \rightarrow \fbox{\rule[-3mm]{0cm}{0.8cm}Data-driven dynamical system} \rightarrow \fbox{\rule[-3mm]{0cm}{0.8cm}Production function}
$$

Thus, we have established the necessary groundwork to reexamine Bowley's law from a mathematical perspective.

\section[Wage share as a differential invariant of a prolonged group action]{Wage share as a differential invariant\\ of a prolonged group action}
\label{s5}

Now that we have established that the Cobb--Douglas function is determined by the exponential growth in the variables $K$, $L$, and $Y$ as functions of $t$ (see \eqref{action}), let us explore the impact of this assumption on wage share without considering the Cobb--Douglas function, as done in Section~\ref{s3}. We aim to demonstrate that, in this scenario, the quantities $s_L$ and $s_K$ remain invariants of the group action \eqref{action}, thus confirming Bowley's law.

Rewriting \eqref{action} in terms of the original variables $K$, $L$, and $Y$, we express them as follows: $\bar{K} = K{\rm e}^{at}$, $\bar{L} = L{\rm e}^{bt}$, $\bar{Y} = Y{\rm e}^{ct}$,
where $a$, $b$, and $c$ are non-negative constants. Next, we derive the corresponding infinitesimal generator and express it using the notations adopted in~\cite{Olver6}:
\begin{equation}\label{inf1}
\mathbf{u}_2 = aK \frac{\partial}{\partial K} + bL \frac{\partial}{\partial L} + cY\frac{\partial}{\partial Y}.
\end{equation}

As the formulas for the quantities $s_L$ and $s_K$ given by \eqref{4} and \eqref{5}, respectively, are defined in terms of the partial derivatives $\partial_L$ and $\partial_K$, it becomes evident that we should produce the first prolongation of the infinitesimal group action \eqref{inf1}. Then, we can employ it to prove that wage and capital shares are indeed invariants.

To achieve this goal, we will utilize the mathematical tools presented in~\cite{Olver6} and~\cite{Saunders89}.

Indeed, let $\big(\mathbb{R}^3,\pi,\mathbb{R}^2\big)$ be a trivial bundle, such that $\pi=pr_1$ and $(K,L,Y)$ be adapted coordinates. Then the corresponding jet bundles are $\big(J^1 \pi, \pi_1, \mathbb{R}^2\big)$ and $\big(J^1\pi, \pi_{1,0}, \mathbb{R}^3\big)$, where (see the commuting diagram below)
\begin{equation}
J^1 \pi=\big\{j^1_p \phi\colon p \in \mathbb{R}^2,\, \phi \in \Gamma_p(\pi) \big\}
\end{equation}
with adapted coordinates $(K,L,Y,Y_K,Y_L).$
\[\begin{tikzcd}
J^1 \pi \arrow{r}{\pi_{1,0}} \arrow[swap]{d}{\pi_1} & \mathbb{R}^3 \arrow{d}{\pi} \\
\mathbb{R}^2 \arrow{r}{\rm id} & \mathbb{R}^2,
\end{tikzcd}
\]
where $\pi_1=\pi \circ \pi_{1,0}.$

The first prolongation of $\mathbf{u}_2$ on $\mathbb{R}^3$ is the vector field ${\rm Pr}^1(\mathbf{u}_2)$ which is a symmetry of the Cartan distribution on $J^1 \pi,$ that is
\begin{align*}
\operatorname{Pr}^1 (\mathbf{u}_2)= aK \frac{\partial}{\partial K}+bL \frac{\partial}{\partial L}+cY\frac{\partial}{\partial Y}+\xi_1(K,L,Y,Y_K,Y_L) \frac{\partial}{\partial Y_K} +\xi_2(K,L,Y,Y_K,Y_L)\frac{\partial}{\partial Y_L}.
\end{align*}

In what follows, the components $\xi_1(K,L,Y,Y_K,Y_L)$ and $\xi_2(K,L,Y,Y_K,Y_L)$ are abbreviated as $\xi_1$ and $\xi_2$ respectively. They are to be determined.

Let us consider a basic contact form $\omega={\rm d}Y-Y_K{\rm d}K-Y_L{\rm d}L$. We require that the one-form~\smash{$\mathcal{L}_{\mathbf{u}^{(1)}_2} (\omega)$} be a contact form \cite{Olver2, Saunders89}. Here $\mathcal{L}$ denotes the Lie derivatives operator. Hence, we compute
\begin{align}
\mathcal{L}_{\mathbf{u}^{(1)}_2} (\omega) = {}& \mathcal{L}_{\mathbf{u}^{(1)}_2} ({\rm d}Y-Y_K{\rm d}K-Y_L{\rm d}L) \nonumber\\
= {}& \mathcal{L}_{\mathbf{u}^{(1)}_2}({\rm d}Y)-\big(\mathcal{L}_{\mathbf{u}^{(1)}_2} Y_K\big){\rm d}K-Y_K \big(\mathcal{L}_{\mathbf{u}^{(1)}_2}({\rm d}K)\big)-\big(\mathcal{L}_{\mathbf{u}^{(1)}_2} Y_L\big){\rm d}L -Y_L \big(\mathcal{L}_{\mathbf{u}^{(1)}_2}({\rm d}L)\big)\nonumber\\
= {}& {\rm d}\big(\mathbf{u}^{(1)}_2(Y)\big)-\big(\mathbf{u}^{(1)}_2(Y_K)\big){\rm d}K-Y_K{\rm d}\big(\mathbf{u}^{(1)}_2(K)\big)-\big(\mathbf{u}^{(1)}_2(Y_L)\big){\rm d}L -Y_L{\rm d}\big(\mathbf{u}^{(1)}_2(L)\big) \nonumber\\
= {}& c {\rm d}Y-\xi_1{\rm d}K-a Y_K{\rm d}K-\xi_2{\rm d}L-b Y_L{\rm d}L\nonumber \\
= {}& c(\omega+Y_K{\rm d}K+Y_L{\rm d}L)-\xi_1{\rm d}K-a Y_K{\rm d}K-\xi_2{\rm d}L-b Y_L{\rm d}L\nonumber\\
= {}& c \omega+(c Y_K-\xi_1-aY_K){\rm d}K+(cY_L-\xi_2-bY_L){\rm d}L,\label{eq1}
\end{align}
the last line of the equation \eqref{eq1} implies that the expressions in the parentheses above vanish, which entails
\[
\xi_1=\xi_1(K,L,Y,Y_K,Y_L)=(c-a)Y_K,\qquad
\xi_2=\xi_2(K,L,Y,Y_K,Y_L)=(c-b)Y_L.
\]

Hence the first prolongation of $\mathbf{u}_2$ is given by
\begin{equation}
\label{first}
{\rm Pr}^{1}(\mathbf{u}_2)=aK \frac{\partial}{\partial K}+bL \frac{\partial}{\partial L}+cY\frac{\partial}{\partial Y}+(c-a)Y_K \frac{\partial}{\partial Y_K} +(c-b)Y_L \frac{\partial}{\partial Y_L}.
\end{equation}
The differential operator \eqref{first} represents the infinitesimal action of a one-parameter Lie group in a 5-dimensional space. According to the fundamental theorem on invariants of regular Lie group actions (see \cite[Chapter~8]{Olver7}), we should expect to determine $5-1 = 4$ fundamental invariants. Indeed, solving the partial differential equation ${\rm Pr}^{1}(\mathbf{u}_2)(F) = 0$ using the method of characteristics, we find
$F = F(I_1, I_2, I_2, I_4)$, where the fundamental invariants $I_i$, $i=1,\ldots, 4$, are found to be
\begin{gather}
\label{i1} I_1=LK^{-\frac{b}{a}},\qquad
I_2=Y K^{-\frac{c}{a}},\qquad
I_3=Y_K K^{\frac{a-c}{a}},\qquad
I_4=Y_LK^{\frac{b-c}{a}}. 
\end{gather}
Next, we combine the invariants \eqref{i1} as appropriate, eliminating the parameters $a$, $b$, and~$c$, to find
\begin{gather}\label{in1}
s_L=\frac{I_1 \cdot I_4}{I_2}=\frac{Y_L L}{Y},
\\
\label{in2}
s_K=\frac{I_3}{I_2}=\frac{Y_K K}{Y}.
\end{gather}
Note that the equation \eqref{in1} provides the classical wage share $s_L$ \eqref{4}, and the equation \eqref{in2} represents the classical capital share $s_K$.

Next, we observe that Sato (see \cite[Chapter 4]{Sato81}) assumed the existence of a homogeneous production function of the form $Y = f(K, L, t)$ and then demonstrated that the wage share defined in terms of the projective variables $x: = L/K$ and $y:= Y/K$, namely $s_L = xy'/y$, was an invariant preserved by an extended action of a vector field generating exponential growth for~$x$ and $y$. The homogeneity of $Y$ allowed for the existence of such variables. For instance, if we have $Y = AL^{\beta}K^{1-\beta}$, it is easy to see that in terms of $x: = L/K$ and $y:= Y/K$, the production function assumes the form $y = Ax^{\beta}$. Under this representation, it becomes an invariant of the infinitesimal action $
U = \frac{1}{1-\beta}x\frac{\partial}{\partial x}+\frac{\beta}{\beta-1}y\frac{\partial}{\partial y}$.
Furthermore, it is possible to show that $s_L = xy'/y$ is also an invariant of the first prolongation of the infinitesimal action defined by $U$. In our previous work \cite{SW20}, we extended this approach to define the notion of modified wage share by assuming logistic growth instead of exponential for $x$ and $y$.

In this paper, we have completely abandoned the need for making any assumptions about a~production function, instead utilizing symmetry methods to compute the first prolongation~\eqref{first} in the original variables $K$, $L$, and $Y$. Indeed, our model does not explicitly require a~production function to be defined. Specifically, we have demonstrated that exponential growth in $K$, $L$, and $Y$ as functions of time \eqref{action5} implies Bowley's law.

\section[The logistic growth model as a generalization of the exponential model]{The logistic growth model as a generalization\\ of the exponential model}
\label{s6}

In Section \ref{s4}, we have established that the Cobb--Douglas production function is a consequence of exponential growth in output and input factors, i.e.,
\[\text{exponential growth} \ \Rightarrow \ \text{the Cobb--Douglas function}.\]
This fact has important ramification for understanding the limitations of the Cobb--Douglas function when it comes to building mathematical and statistical models in economic growth theory. Viewing the Cobb--Douglas function as a consequence of a vigorous growth in production, capital, and labor enables us, for example, to understand better why so many data sets describing the current state of many economies cannot be accurately described by this particular production function (see, for example, Antr\`{a}s~\cite{PA04}, Xiang~\cite{HX04}, and Gechert et al.~\cite{GHIK19}). It is also evident that labor no longer grows exponentially (see, for example, Aghion and Howitt~\cite{AH08}).

We have extended in \cite{SW20} the exponential model by replacing the assumption about exponential growth in labor, capital, and production with the corresponding assumption that labor, capital, and production grow logistically, arriving, as a result, at the following dynamical system:%
\begin{gather}\label{model2}
\dot{x}_i = b_ix_i\left(1 - \frac{x_i}{N_i}\right), \qquad i = 1, 2, 3,
\end{gather}
where $x_1 = L$ (labor), $x_2 = K$ (capital), $x_3 = Y$ (production) and the constants $N_i$, ${i=1,2,3}$, denote the corresponding carrying capacities. Furthermore, we employed Sato's approach~\cite{Sato81} to integrate this dynamical system and thus derive a new production function (see \cite[formula~(4.5)]{SW20}).

Since the dynamical system \eqref{model2} reduces to \eqref{model1} as $N_i \to\infty$, $i=1,2,3$, we can naturally consider the exponential model as a limiting case of the logistic one, which also means that the Cobb--Douglas function \eqref{CD} can be viewed as a limiting case of the new production function derived in \cite{SW20} from the logistical model.

Indeed, integrating \eqref{model2}, we arrive at the corresponding one-parameter Lie group action given by
\begin{equation}
\label{action2}
x_i = \frac{N_i x_i^0}{x_i^0 + \big(N_i - x^0_i\big){\rm e}^{-b_it}}, \qquad i = 1,2,3,
\end{equation}
where $x_i^0$, $i = 1,2,3$ are the initial conditions. It follows from \eqref{action2}
\[\frac{x_i \big(N_i - x_i^0\big)}{x_i^0(N_i - x_i)} = {\rm e}^{b_it}, \qquad i=1,2,3.\]
Next, we obtain
\begin{equation}
\label{newfunction3}
\prod_{i=1}^3\left[\frac{x_i \big(N_i - x_i^0\big)}{x_i^0(N_i - x_i)}\right]^{a_i} = {\rm e}^{\left(\sum_{i=1}^3\alpha_ib_i\right)t},
\end{equation}
where $a_i$, $i=1,2,3$ are some parameters. We see that the left-hand side of the equation~\eqref{newfunction3} is an invariant of the one-parameter group action generated by \eqref{model2} if and only if the parameters $a_i$, $i=1,2,3$ satisfy the same orthogonality condition~\eqref{linear} for the fixed values of the parameters~$b_i$, $i=1,2,3$ determined by~\eqref{model2}. In terms of the language of induced mappings the natural connection between the two models appears even simpler. Indeed, introduce the following map
 $\psi\colon \mathbb{R}_+^3 \rightarrow D,$ $ D = ]0,N_1[{}\times{}]0,N_2[{}\times{}]0,N_3[ {}\subset \mathbb{R}^3$ via the coordinate transformation
\begin{equation}
\label{newcoordinates}
\tilde{x}_i = \frac{N_i x_i}{N_i + x_i}, \qquad i = 1,2,3.
\end{equation}
We note that $\psi$ is a diffeomorphism, $\tilde{x}_i \to x_i$, as $N_i \to \infty$, $0<\tilde{x}_i < N_i$, $i=1,2,3$, and the corresponding Jacobian is given by
\begin{equation}
\label{J}
J_{\psi}(x,\tilde{x}) = \frac{\partial \tilde{x}_i}{\partial x_j} = \mbox{diag} \left(\frac{N_1^2}{(N_1+x_1)^2}, \frac{N_2^2}{(N_2+x_2)^2}, \frac{N_3^2}{(N_3+x_3)^2}\right), \qquad i,j=1,2,3.
\end{equation}

In this view, the logistic model defined above reduces to previous exponential model as $N_i \to\infty$, $i=1,2,3$. Indeed, in terms of the new coordinates $\tilde{x}_i$, $i=1,2,3$ given by \eqref{newcoordinates} the formula \eqref{newfunction3} reduces to \eqref{GCDF}. Alternatively, the logistic model is simply the image of the exponential model under the push-forward map $\psi_*\colon T_{x}\mathbb{R}^3_+ \rightarrow T_{\psi (x)}\mathbb{R}^3_+$ induced by $\psi$. Clearly, the vector field
\[X = \sum_{i=1}^3b_ix_i \frac{\partial}{\partial x_i}\]
that defines the system \eqref{model1} is mapped by $\psi_*$ to the vector field $\tilde{X} = \tilde{X}^i\frac{\partial }{\partial \tilde{x}_i}$, where
\begin{equation}
\label{logistic}
\tilde{X}^i = X^j \frac{\partial \tilde{x}_i}{\partial x_j} = b_i\frac{N_i\tilde{x}_i}{N_i-\tilde{x}_i}\frac{N_i^2}{\big(N_i+\frac{N_i\tilde{x}_i}{N_i-\tilde{x}_i}\big)^2}=b_i\tilde{x}_i\left(1-\frac{\tilde{x}_i}{N_i}\right),
\end{equation}
where we have used the formulas \eqref{newcoordinates} and \eqref{J}. Clearly, the vector field $\tilde{X}$ defined by the components \eqref{logistic} is precisely the vector field that gives rise to the system \eqref{model2}, as expected. Similarly, a smooth function $f\colon\mathbb{R}^3_+\to \mathbb{R}$ is pulled back to a smooth function $\psi^*f$ defined by~${(\psi^* f)(x) = f(\psi (x))}$. Therefore, substituting \smash{${x}_i = \frac{N_i\tilde{x}_i}{N_i-\tilde{x}_i}$} in \eqref{CDF1}, dropping tildes, and then solving for $x_3$, we arrive at the following function $f\colon {}]0,N_1[{}\times{} ]0,N_2[ {}\rightarrow{} ]0,N_3[$ given by
\begin{equation}
\label{newfunction1}
x_3 = f(x_1, x_2) = \frac{N_3x_1^{\alpha}x_2^{\beta}}{\frac{N_3}{N_1^{\alpha}N_2^{\beta}B}(N_1-x_1)^{\alpha}(N_2-x_2)^{\beta} + x_1^{\alpha}x_2^{\beta}},
\end{equation}
where
\[B = A\left(\frac{N_1-x_1^0}{N_1x_1^0}\right)^{\frac{a_1}{a_3}}\left(\frac{N_2-x_2^0}{N_2x_2^0}\right)^{\frac{a_2}{a_3}}\left(\frac{N_3-x_3^0}{N_3x_3^0}\right),
\]
$\alpha = -\frac{a_1}{a_3}$, and $\beta = -\frac{a_2}{a_3}$.
It is easy to show that the Cobb--Douglas function \eqref{CDF1} is the limiting case of the function \eqref{newfunction1}, as $N_i \to \infty$, $i=1,2,3$.

 Now, identifying $x_3 = Y$ (production), $x_1 = L$ (labor), $x_2 = K$ (capital), $N_1 = N_L$, $N_2 = N_K$, $N_3 = N_Y$, \smash{$C = N_3N_1^{-\alpha}N_2^{-\beta}B^{-1}$}, and extending the domain and range of the function~\eqref{newfunction1} to~$\mathbb{R}_+^2$ and~$\mathbb{R}_+$ respectively, we arrive at the production function
\begin{equation}
\label{newfunction2}
Y = f(L, K) = \frac{N_YL^{\alpha}K^{\beta}}{C|N_L-L|^{\alpha}|N_K- K|^{\beta} + L^{\alpha}K^{\beta}}
\end{equation}
derived for the first time in \cite{SW20} by employing integration.

We note that the one-input version of the production function \eqref{newfunction2}, namely the function
\[
Y = f(x) = \frac{N_fx^{\alpha}}{C|N_x-x|^{\alpha} + x^{\alpha}},
\]
derived in \cite{SW20}, apparently gives a mathematical formula defining the environmental Kuznets curve (see, for example, \cite[Figure~1]{DIS04}).

By analogy with the generalization of the Cobb--Douglas function presented in Section~\ref{s4}, we propose Table~\ref{table2}.

\begin{table}\renewcommand{\arraystretch}{1.1}\centering
\begin{tabular}{ |l||c|c|c|c| }
\hline
& $b$ & $x^0$ & $N$ & RSS\\
\hline\hline
labor &0.07842367 &2.092004 &175.97 &508.0948\\ [0.1cm]
capital & 0.07793777 & 1.575667 &230.26 &299.7033\\ [0.1cm]
production & 0.04619786 &11.312991 &211.30 &419.7767\\
\hline
\end{tabular}
\caption{The estimated values of the parameters that determine logistic growth \eqref{action2}.}\label{table2}
\end{table}

\begin{Definition}[logistic production function] \label{d2} Given the one-parameter group action
\begin{equation}\label{action7}
x_i = \frac{N_i x_i^0}{x_i^0 + \big(N_i - x^0_i\big){\rm e}^{-b_it}}, \qquad x_i^0, b_i, N_i >0, \qquad i = 1, \ldots, n
\end{equation}
in $\mathbb{R}^n_+$. Then the {\em logistic production function} is defined as an element of the following family of invariants of the one-parameter Lie group action \eqref{action7}:
\begin{equation}
\label{LPF}
\prod_{i=1}^n \left[\frac{x_i \big(N_i - x_i^0\big)}{x_i^0(N_i - x_i)}\right]^{a_i} = C, \qquad a_i, N_i > 0 \qquad i = 1, \ldots, n,
\end{equation}
where $C \in \mathbb{R}$ is an arbitrary constant and $x_i^0$, $i = 1, \ldots, n$ are the corresponding initial conditions, provided the orthogonality condition
\begin{equation}
\label{linear6}
\sum_{i=1}^n a_ib_i = 0
\end{equation}
holds true.
\end{Definition}

Therefore, we have shown that the two models, namely, the exponential model, discussed in the preceding section, and the logistic model are, in fact, equivalent, modulo the transformation~$\psi$ given by \eqref{newcoordinates}.

Now let us revisit again the data studied in \cite{SW20}, i.e., the data from the period 1947--2016 provided by the Federal Reserve Bank of St. Louis (\url{https://fred.stlouisfed.org}), employing the FRED tool. Recall that in~\cite{SW20} we fitted the function \eqref{newfunction2} to the data, assuming the values of the carrying capacities to be $N_Y = 120$, $N_L = N_K = 150$, and the condition~\eqref{7}. In what follows, we employ a more delicate and accurate approach. First, we note that the condition~\eqref{7} is redundant in this case, since even when $\alpha + \beta =1$ the production function~\eqref{newfunction2} is not homogeneous. Second, dropping the condition~\eqref{7},
now we start by fitting the formulas of logistic growth~\eqref{action2} to the data representing production, labor and capital, rather than the formula~\eqref{newfunction2} itself. More specifically, we treat the parameters $b_i$, $x^0_i$, $N_i$, $i=1,2,3$ in~\eqref{action2} as the predictors of logistic growth and, with the aid of R and the least squares method, arrive at the following values, presented along with the corresponding RSS's in Table~\ref{table2}.

\begin{figure}[t]\centering
\includegraphics[width=0.53\textheight]{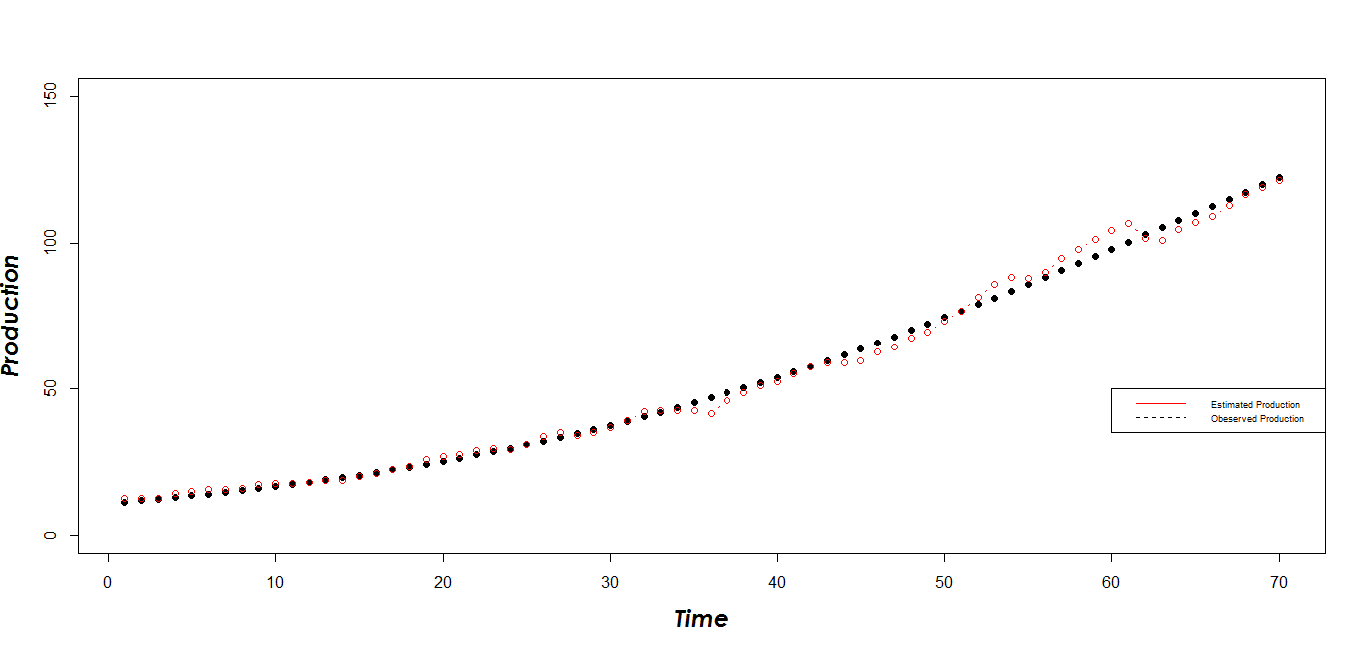}
\caption{Observed vs estimated production.}\label{figure12}
\end{figure}

\begin{figure}[ht]
\centering
\includegraphics[width=0.53\textheight]{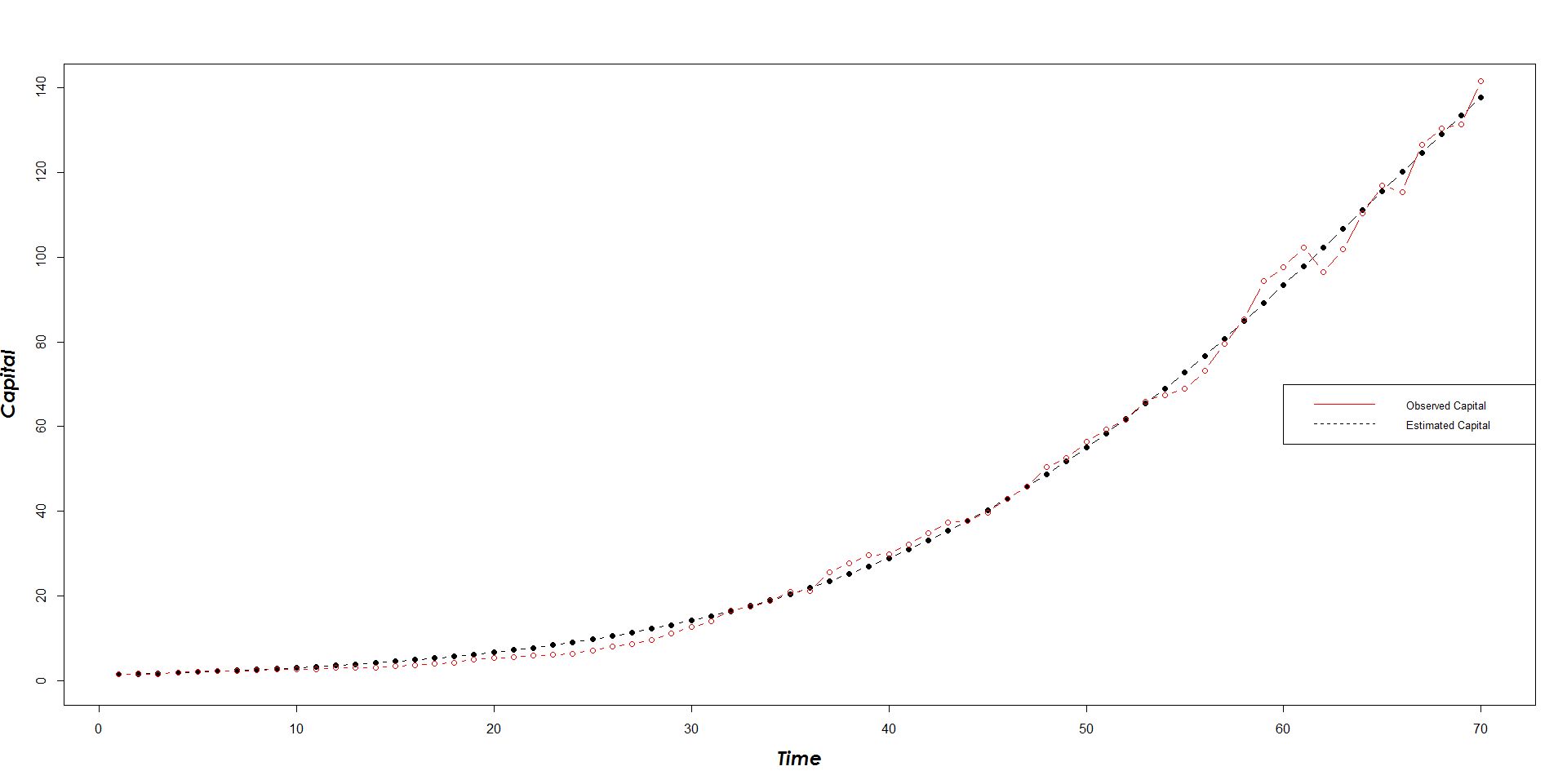}
\caption{Observed vs estimated capital.}\label{figure13}
\end{figure}

\begin{figure}[ht]\centering
\includegraphics[width=0.53\textheight]{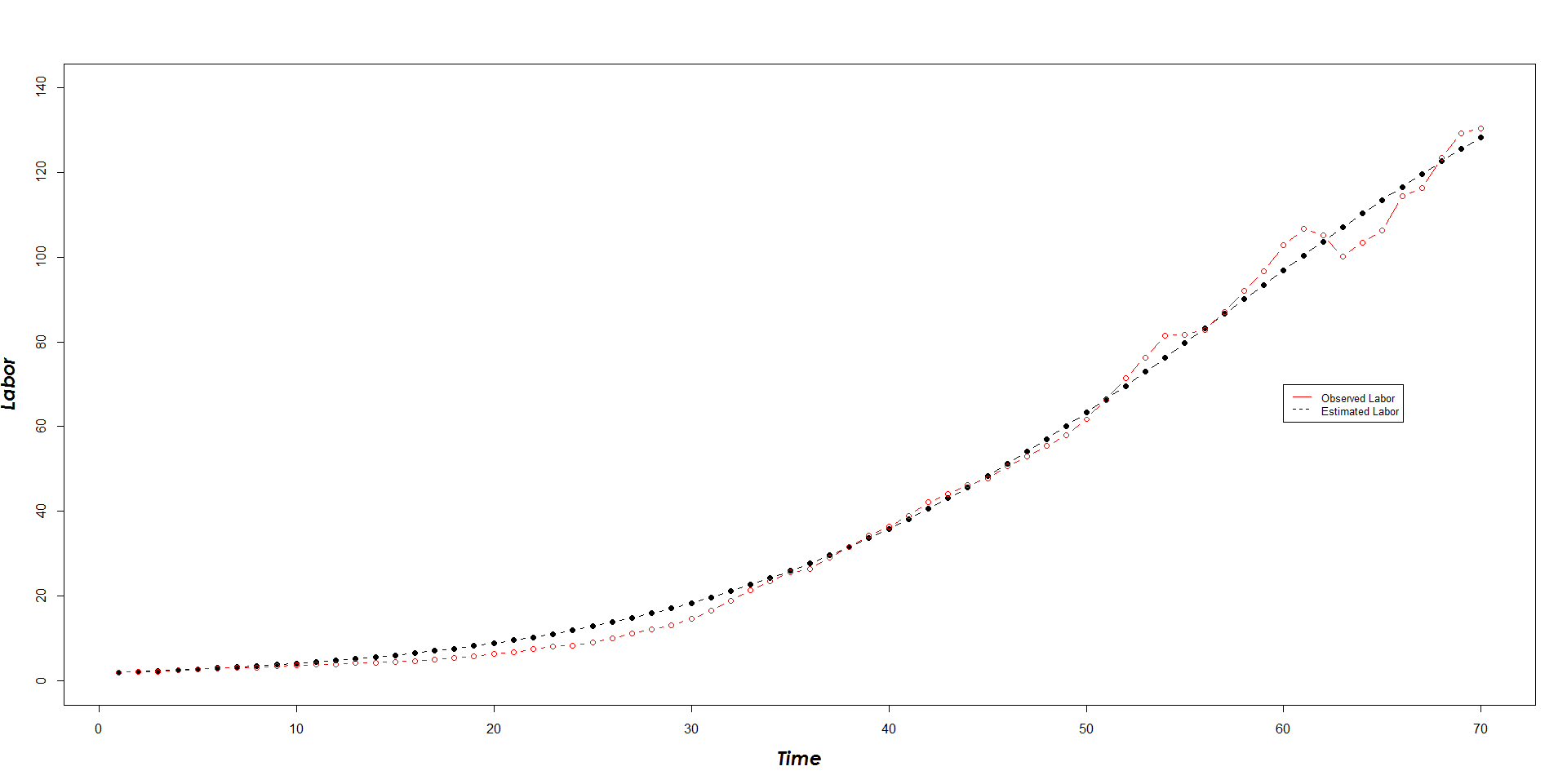}
\caption{Observed vs estimated labor.}
\label{figure14}
\end{figure}

We first note that the RSS values tabulated in Table \ref{table2} are much better than the ones obtained for the fixed values of the carrying capacities assumed in \cite{SW20}. We illustrate the accuracy of our estimations by Figures~\ref{figure12},~\ref{figure13} and~\ref{figure14}.

 Next, using the carrying capacities presented in Table~\ref{table2} and dropping the condition~\eqref{7}, we arrive at the following estimated values for the parameters of the production function~\eqref{newfunction2}:
\begin{equation}
\label{alphabeta3}
\alpha=0.46780229, \qquad \beta=0.05955408, \qquad C=1.59899336
\end{equation}
that make it fit to the data with ${\rm RSS}=428.27$. We note that the parameters $b_1$~(labor), $b_2$~(capital), $b_3$ (production) presented in Table \ref{table2} and the parameters $\alpha$, $\beta$ in \eqref{alphabeta3} satisfy the orthogonality condition \eqref{linear} with good accuracy, as expected. It must be noted at this point that the function \eqref{newfunction2} with $N_Y = 120$ and $N_L = N_K = 150$ under the condition \eqref{7} yielded a much less accurate fit to this data in \cite{SW20} with ${\rm RSS} = 4336.976$. Finally, we compute the corresponding RSS for the Cobb--Douglas function \eqref{6} determined by the exponential growth parameters \eqref{data1} fitted to the same data from the period 1947--2016 provided by the Federal Reserve Bank of St.~Louis (\url{https://fred.stlouisfed.org}), employing the FRED tool, to find that in this case RSS = 584.4616 which is less than the RSS for the function~\eqref{newfunction2} determined by the parameters presented in Table~\ref{table2}.

Recall that a similar, ``S-shaped" production function
\begin{gather}
\label{Sshaped}
Y = g(K, L) = \frac{a K^pL^{1-p}}{1 + bK^pL^{-p}}
\end{gather}
was recently introduced, employing a heuristic approach, see, for example, Anita et al.~\cite{ACKL13}, Capasso et al.\ \cite{CEL12}, Engbers et al.~\cite{EBC14}, La Torre et al.~\cite{LLM15} and the relevant references therein for more details and applications. Note that the production function \eqref{Sshaped} is reducible to the Cobb--Douglas function~\eqref{CD} (i.e., when $b = 0$). Compare the function \eqref{Sshaped} to the Holling disc equation in eco-dynamics \cite{CSH59}. Also, we observe that the new production function~\eqref{newfunction2} is reducible to the production function \eqref{Sshaped} when $K$ and $L$ $\ll$ $N_K$ and $N_L$ respectively, $N_L, N_K \approx 1$, $C = 1$ in~\eqref{newfunction2} and $a = N_{Y}$, $b = 1$ in \eqref{Sshaped}.

\section{The labor share for the logistic model}\label{s7}

Now let us investigate what transpires when the underlying model is logistic, rather than exponential. Indeed, substituting the formula for the logistic production function~\eqref{newfunction2} into~\eqref{4}, we arrive at the following expression for the labor (wage) share:
\begin{equation}
\label{laborshare1}
s_L = \alpha \frac{N_L}{N_L-L}\cdot \frac{C|N_L-L|^{\alpha}|N_K-K|^{\beta}}{C|N_L-L|^{\alpha}|N_K-K|^{\beta}+L^{\alpha}K^{\beta}}.
\end{equation}
Note that as $N_L, N_K \to \infty$, the labor share $s_L \to \alpha$, where $s_L$ is given by \eqref{laborshare1}, as expected. In view of the formulas \eqref{action2}, $s_L$ assumes the following form in terms of the parameter $t$:
\begin{equation}
\label{laborshare2}
s_L=\frac{b_3}{b_1} \cdot \frac{(N_Y-Y_0){\rm e}^{(b_1-b_3) t}}{N_L-L_0} \cdot \frac{L_0+(N_L-L_0){\rm e}^{-b_1t}}{Y_0+(N_Y-Y_0){\rm e}^{-b_3t}}.
\end{equation}
The above formulas \eqref{laborshare1} and \eqref{laborshare2} put in evidence that the labor share $s_L$ is no longer constant for the logistic model, which illustrates the limitations of Bowley's law.

\section{Concluding remarks}
\label{s8}

For many years, Bowley's law has stood as a cornerstone in economics, often observed through the analysis of data from specific countries or regions. Indeed, modern growth and distribution theories have been founded on the assumption of Bowley's law holding true.
However, recent times have seen a surge in critics challenging Bowley's law. Many researchers, scrutinizing data from the past few decades, have noted a decline in the wage share over time, sparking significant controversy among economists. Some argue that Bowley's law might never have been valid in the first place (see Kr\"{a}mer \cite{Kramer11}, for example).

In our paper, we developed a model that combines statistical data and extends the analytical model proposed by Sato \cite{Sato81}. Our model demonstrates that Bowley's law primarily arises as a~consequence of exponential growth in capital ($K$), labor ($L$), and production ($Y$) as functions of time $t$.

Importantly, our model does not rely on assumptions about the specific parameters governing the exponential growth of $K$, $L$, and $Y$, or the form of the associated Cobb--Douglas-type production functions.
Thus, it becomes evident that Bowley's law is indeed a consequence of exponential growth and nothing else. In other words, for Bowley's law to hold true, the underlying economy must exhibit rigorous growth.

On the other hand, if Bowley's law does not hold true, it indicates that $K$, $L$, and $Y$ are not all growing exponentially. We demonstrated this fact by employing the logistic model and the corresponding production function. For the logistic model, for example, the labor share $s_L$ is no longer an invariant, as expected (see \eqref{laborshare2}). We have also shown that the exponential model can be viewed as a limiting case of the exponential model.

In addressing this issue, we employed mathematical techniques, particularly symmetry methods, which shed some light on this problem and contributed to resolving, at least to some extent, the controversy surrounding this economic law.

\subsection*{Acknowledgements}

The authors sincerely appreciate the invaluable mathematical knowledge they have acquired from Peter Olver, which has been instrumental in the development of this and other projects in the past. The first author (RGS) would like to offer special thanks to the organizers of the conference ``Symmetry, Invariants, and their Applications: A Celebration of Peter Olver's 70th Birthday'' for extending an invitation to participate and present a talk. The conference provided an exceptional platform for engaging in stimulating discussions with fellow participants.
Additionally, the authors express their gratitude to the anonymous referees for their valuable feedback, new references, insightful suggestions, and constructive critiques, all of which have significantly enhanced the quality of the presentation.

\pdfbookmark[1]{References}{ref}
\LastPageEnding


\begin{thebibliography}{99}
\footnotesize\itemsep=-0.5pt

\bibitem{AH08}
Aghion P., Howitt P., The economics of growth, The MIT Press, Cambridge, MA,
 2008.

\bibitem{ACKL13}
Ani\c{t}a S., Capasso V., Kunze H., La~Torre D., Optimal control and long-run
 dynamics for a spatial economic growth model with physical capital
 accumulation and pollution diffusion, \href{https://doi.org/10.1016/j.aml.2013.04.002}{\textit{Appl. Math. Lett.}} \textbf{26}
 (2013), 908--912.

\bibitem{PA04}
Antr\`as P., Is the {US} aggregate production function {C}obb--{D}ouglas? {N}ew
 estimates of the elasticity of substitution, \href{https://doi.org/10.2202/1534-6005.1161}{\textit{Contrib. Macroecon.}}
 \textbf{4} (2004), 1--34.

\bibitem{Beaudreau17}
Beaudreau B.C., The economies of speed {$KE = 1/2mv^2$} and the productivity
 slowdown, \href{https://doi.org/10.1016/j.energy.2017.02.022}{\textit{Energy}} \textbf{124} (2017), 100--113.

\bibitem{Bowley1900}
Bowley A.L., Wages in the {U}nited {K}ingdom in the nineteenth century: {N}otes
 for the use of students of social and economic questions, Cambridge
 University Press, Cambridge, 1900.

\bibitem{Bowley1937}
Bowley A.L., Wages and income in the {U}nited {K}ingdom since 1860, Cambridge
 University Press, Cambridge, 1937.

\bibitem{CEL12}
Capasso V., Engbers R., La~Torre D., Population dynamics in a spatial {S}olow
 model with a convex-concave production function, in Mathematical and
 {S}tatistical {M}ethods for {A}ctuarial {S}ciences and {F}inance, \href{https://doi.org/10.1007/978-88-470-2342-0_8}{Springer},
 Dordrecht, 2012, 61--68.

\bibitem{CS21}
Cherevatskyi D., Smirnov R.G., A novel approach to characterizing the
 relationship between economic growth and energy consumption, \href{https://doi.org/10.15407/economyukr.2021.12.057}{\textit{Econ.
 Ukraine}} (2021), no.~12, 57--70.

\bibitem{CD28}
Cobb C.W., Douglas P.H., A theory of production, \textit{Amer. Econ. Rev.}
 \textbf{18} (1928), 139--165.

\bibitem{CMS17}
Cochran C.M., McLenaghan R.G., Smirnov R.G., Equivalence problem for the
 orthogonal separable webs in 3-dimensional hyperbolic space, \href{https://doi.org/10.1063/1.4983998}{\textit{J.~Math.
 Phys.}} \textbf{58} (2017), 063513, 43~pages.

\bibitem{Douglas76}
Douglas P.H., The {C}obb--{D}ouglas production function once again: {I}ts
 history, its testing, and some new empirical values, \href{https://doi.org/10.1086/260489}{\textit{J.~Polit. Econ.}}
 \textbf{84} (1976), 903--916.

\bibitem{EBC14}
Engbers R., Burger M., Capasso V., Inverse problems in geographical economics:
 parameter identification in the spatial {S}olow model, \href{https://doi.org/10.1098/rsta.2013.0402}{\textit{Philos.
 Trans.~R.~Soc. Lond. Ser.~A Math. Phys. Eng. Sci.}} \textbf{372} (2014),
 20130402, 13~pages.

\bibitem{Fukang96}
Fang F., The symmetry approach on economic systems, \href{https://doi.org/10.1016/S0960-0779(96)00083-5}{\textit{Chaos Solitons
 Fractals}} \textbf{7} (1996), 2247--2257.

\bibitem{Olver6}
Fels M., Olver P.J., Moving coframes.~{I}. {A}~practical algorithm,
 \href{https://doi.org/10.1023/A:1005878210297}{\textit{Acta Appl. Math.}} \textbf{51} (1998), 161--213.

\bibitem{Olver6a}
Fels M., Olver P.J., Moving coframes.~{II}. {R}egularization and theoretical
 foundations, \href{https://doi.org/10.1023/A:1006195823000}{\textit{Acta Appl. Math.}} \textbf{55} (1999), 127--208.

\bibitem{Olver4}
Fokas A.S., Olver P.J., Rosenau P., A plethora of integrable bi-{H}amiltonian
 equations, in Algebraic {A}spects of {I}ntegrable {S}ystems, \textit{Progr.
 Nonlinear Differential Equations Appl.}, Vol.~26, \href{https://doi.org/10.1007/978-1-4612-2434-1_5}{Birkh\"auser}, Boston, MA,
 1997, 93--101.

\bibitem{GHIK19}
Gechert S., Havranek T., Irsova Z., Kolcunova D., Death to the
 {C}obb--{D}ouglas production function? {A}~quantitative survey of the
 capital-labor substitution elasticity, {E}conStor {P}reprints, {ZBW},
 {L}eibniz {I}nformation {C}entre for {E}conomics, 2009.

\bibitem{CSH59}
Holling C.S., Some characteristics of simple types of predation and parasitism,
 in The {C}anadian {E}ntomologist, \href{https://doi.org/10.4039/Ent91385-7}{Cambridge University Press}, Canada, 1959,
 385--398.

\bibitem{HMS05}
Horwood J.T., McLenaghan R.G., Smirnov R.G., Invariant classification of
 orthogonally separable {H}amiltonian systems in {E}uclidean space,
 \href{https://doi.org/10.1007/s00220-005-1331-8}{\textit{Comm. Math. Phys.}} \textbf{259} (2005), 679--709,
 \href{https://arxiv.org/abs/math-ph/0605023}{arXiv:math-ph/0605023}.

\bibitem{Humphrey97}
Humphrey T.M., Algebraic production functions and their uses before
 {C}obb--{D}ouglas, \textit{FRB Richmond Econ. Quart.} \textbf{83} (1997),
 51--83.

\bibitem{Olver11}
Kamran N., Olver P.J., Tenenblat K., Local symplectic invariants for curves,
 \href{https://doi.org/10.1142/S0219199709003326}{\textit{Commun. Contemp. Math.}} \textbf{11} (2009), 165--183.

\bibitem{Olver13}
Kogan I.A., Olver P.J., Invariants of objects and their images under surjective
 maps, \href{https://doi.org/10.1134/S1995080215030063}{\textit{Lobachevskii~J. Math.}} \textbf{36} (2015), 260--285,
 \href{https://arxiv.org/abs/1509.06690}{arXiv:1509.06690}.

\bibitem{Kramer11}
Kr\"amer H.M., Bowley's law: {T}he diffusion of an empirical supposition into
 economic theory, \href{https://doi.org/10.3917/cep.061.0019}{\textit{Pap. Polit. Econ.}} \textbf{61} (2011), 19--49.

\bibitem{KKTT19}
Kreinovich V., Kosheleva O., Limit theorems as blessing of dimensionality:
 neural-oriented overview, \href{https://doi.org/10.3390/e23050501}{\textit{Entropy}} \textbf{23} (2021), 501, 19~pages.

\bibitem{LLM15}
La~Torre D., Liuzzi D., Marsiglio S., Pollution diffusion and abatement
 activities across space and over time, \href{https://doi.org/10.1016/j.mathsocsci.2015.09.001}{\textit{Math. Social Sci.}} \textbf{78}
 (2015), 48--63.

\bibitem{Olver1}
Olver P.J., Evolution equations possessing infinitely many symmetries,
 \href{https://doi.org/10.1063/1.523393}{\textit{J.~Math. Phys.}} \textbf{18} (1977), 1212--1215.

\bibitem{Olver2}
Olver P.J., Applications of {L}ie groups to differential equations,
 \textit{Grad. Texts in Math.}, Vol. 107, \href{https://doi.org/10.1007/978-1-4684-0274-2}{Springer}, New York, 1986.

\bibitem{Olver3}
Olver P.J., Equivalence, invariants, and symmetry, \href{https://doi.org/10.1017/CBO9780511609565}{Cambridge University Press},
 Cambridge, 1995.

\bibitem{Olver7}
Olver P.J., Classical invariant theory, \textit{London Math. Soc. Stud. Texts},
 Vol.~44, \href{https://doi.org/10.1017/CBO9780511623660}{Cambridge University Press}, Cambridge, 1999.

\bibitem{Olver8}
Olver P.J., Moving frames~-- in geometry, algebra, computer vision, and
 numerical analysis, in Foundations of {C}omputational {M}athematics
 ({O}xford, 1999), \textit{London Math. Soc. Lecture Note Ser.}, Vol. 284,
 Cambridge University Press, Cambridge, 2001, 267--297.

\bibitem{Olver12}
Olver P.J., Introduction to partial differential equations, Undergrad. Texts
 Math., \href{https://doi.org/10.1007/978-3-319-02099-0}{Springer}, Cham, 2014.

\bibitem{Olver15}
Olver P.J., Motion and continuity, \href{https://doi.org/10.1007/s00283-022-10194-x}{\textit{Math. Intelligencer}} \textbf{44}
 (2022), 241--249.

\bibitem{Olver10}
Olver P.J., Pohjanpelto J., Moving frames for {L}ie pseudo-groups,
 \href{https://doi.org/10.4153/CJM-2008-057-0}{\textit{Canad.~J. Math.}} \textbf{60} (2008), 1336--1386.

\bibitem{Olver9}
Olver P.J., Shakiban C., Applied linear algebra, Pearson Prentice Hall, Upper
 Saddle River, NJ, 2006.

\bibitem{Olver14}
Olver P.J., Valiquette F., Recursive moving frames for {L}ie pseudo-groups,
 \href{https://doi.org/10.1007/s00025-018-0818-5}{\textit{Results Math.}} \textbf{73} (2018), 57, 64~pages.

\bibitem{OG23}
Orlando G., On the assumptions of the {C}obb--{D}ouglas production function and
 their assessment in contemporary economic theory, {a}vailable at
 \url{https://ssrn.com/abstract=4439917}, 2023, 21~pages.

\bibitem{PY18}
Perets G., Yashiv E., Lie symmetries and essential restrictions in economic
 optimization, {CEPR} {D}iscussion {P}aper {N}o.~{DP}12611, 2018, 33~pages.

\bibitem{PS05}
Praught J., Smirnov R.G., Andrew {L}enard: a mystery unraveled, \href{https://doi.org/10.3842/SIGMA.2005.005}{\textit{SIGMA}}
 \textbf{1} (2005), 005, 7~pages, \href{https://arxiv.org/abs/nlin.SI/0510055}{arXiv:nlin.SI/0510055}.

\bibitem{PS1964}
Samuelson P., Economics: {A}n introductory textbook, McGraw-Hill, New York,
 1964.

\bibitem{Sato81}
Sato R., Theory of technical change and economic invariance. {A}pplication of
 Lie groups, Economic Theory, Econometrics, and Mathematical Economics,
 Academic Press, Inc., New York, 1981.

\bibitem{SR14}
Sato R., Ramachandran R.V., Symmetry and economic invariance, \textit{Adv. Jpn.
 Bus. Econ.}, Vol.~1, \href{https://doi.org/10.1007/978-4-431-54430-2}{Springer}, Tokyo, 2014.

\bibitem{Saunders89}
Saunders D.J., The geometry of jet bundles, \textit{London Math. Soc. Lecture
 Note Ser.}, Vol. 142, \href{https://doi.org/10.1017/CBO9780511526411}{Cambridge University Press}, Cambridge, 1989.

\bibitem{Schneider11}
Schneider D., The labor share: {A} review of theory and evidence, {SFB}649
 {E}conomic {R}isk, {D}iscussion {P}aper, 2011.

\bibitem{SW20}
Smirnov R.G., Wang K., In search of a new economic model determined by logistic
 growth, \href{https://doi.org/10.1017/s0956792519000081}{\textit{European~J. Appl. Math.}} \textbf{31} (2020), 339--368,
 \href{https://arxiv.org/abs/1711.02625}{arXiv:1711.02625}.

\bibitem{SW21}
Smirnov R.G., Wang K., The {C}obb--{D}ouglas production function revisited, in
 Proceedings of the International Conference on Applied Mathematics, Modeling
 and Computational Science ``AMMCS-2019'' (August 18--23, 2019, Waterloo,
 Ontario, Canada), \textit{Springer Proc. Math. Stat.}, Vol. 343, \href{https://doi.org/10.1007/978-3-030-63591-6_66}{Springer},
 Cham, 2021, 725--734, \href{https://arxiv.org/abs/1910.06739}{arXiv:1910.06739}.

\bibitem{SWW22}
Smirnov R.G., Wang K., Wang Z., The {C}obb--{D}ouglas production function for
 an exponential model, in Advances in Econometrics, Operational Research, Data
 Science and Actuarial Studies: Techniques and Theories, \href{https://doi.org/10.1007/978-3-030-85254-2_1}{Springer}, Cham, 2022, 1--12.

\bibitem{DIS04}
Stern D.I., The rise and fall of the environmental {K}uznets curve,
 \href{https://doi.org/10.1016/j.worlddev.2004.03.004}{\textit{World Develop.}} \textbf{28} (2004), 1419--1439.

\bibitem{SE13}
Stockhammer E., Why have wage shares fallen? {A}n analysis of the determinants
 of functional income dis\-tri\-bution, in Wage-Led Growth, \textit{Adv. Labour Stud.},
 \href{https://doi.org/10.1057/9781137357939_3}{Palgrave Macmillan}, London, 2013, 40--70.

\bibitem{HX04}
Xiang H., Is {C}anada's aggregate production function {C}obb--{D}ouglas?
 {E}stimation of the elasticity of substitution between capital and labor,
 {M}aster Thesis, {S}imon {F}raser {U}niversity, 2004.

\bibitem{Olver5}
Yezzi A., Kichenassamy S., Kumar A., Olver P.J., Tannenbaum A., A geometric
 snake model for segmentation of medical imagery, \href{https://doi.org/10.1109/42.563665}{\textit{IEEE Trans. Medical
 Imaging}} \textbf{16} (1997), 199--209.

\end{thebibliography}
\end{document}